\documentclass[11pt]{amsart}
\usepackage{mathrsfs}
\usepackage{amssymb}

\usepackage{amscd}
\usepackage{amsmath, amssymb}
\usepackage{amsfonts}
\usepackage[colorlinks,linkcolor=blue,citecolor=blue, pdfstartview=FitH]{hyperref}

\numberwithin{equation}{section}

\theoremstyle{plain}
\newtheorem{thm}{Theorem}[section]
\newtheorem{theorem}[thm]{Theorem}
\newtheorem{lemma}[thm]{Lemma}
\newtheorem{corollary}[thm]{Corollary}
\newtheorem{proposition}[thm]{Proposition}

\theoremstyle{definition}

\newtheorem{remark}[thm]{Remark}

\newtheorem{definition}[thm]{Definition}

\newtheorem{example}[thm]{Example}

\newtheorem{defn-thm}[thm]{Definition-Theorem}

\newcommand{\sE}{{\mathcal E}}

\newcommand{\sL}{{\mathcal L}}

\newcommand{\sO}{{\mathcal O}}

\newcommand{\C}{{\mathbb C}}

\renewcommand{\P}{{\mathbb P}}
\newcommand{\Q}{{\mathbb Q}}
\newcommand{\R}{{\mathbb R}}

\renewcommand{\S}{{\mathbb S}}

\newcommand{\qtq}[1]{\quad\mbox{#1}\quad}
\newcommand{\bp}{\bar{\partial}}
\newcommand{\Om}{\Omega}

\newcommand{\ds}{\oplus}

\newcommand{\ts}{\otimes}

\newcommand{\btheorem}{\begin{theorem}}
\newcommand{\etheorem}{\end{theorem}}
\newcommand{\bproposition}{\begin{proposition}}
\newcommand{\eproposition}{\end{proposition}}
\newcommand{\bdefinition}{\begin{definition}}
\newcommand{\edefinition}{\end{definition}}
\newcommand{\bcorollary}{\begin{corollary}}
\newcommand{\ecorollary}{\end{corollary}}
\newcommand{\bproof}{\begin{proof}}
\newcommand{\eproof}{\end{proof}}
\newcommand{\bremark}{\begin{remark}}
\newcommand{\eremark}{\end{remark}}
\newcommand{\eexample}{\end{example}}
\newcommand{\bexample}{\begin{example}}
\newcommand{\la}{\langle}
\newcommand{\elemma}{\end{lemma}}
\newcommand{\blemma}{\begin{lemma}}
\newcommand{\ra}{\rangle}
\newcommand{\sq}{\sqrt{-1}}

\newcommand{\p}{\partial}

\renewcommand{\bar}{\overline}
\newcommand{\eps}{\varepsilon}

\renewcommand{\phi}{\varphi}

\newcommand{\ee}{\end{eqnarray*}}
\newcommand{\be}{\begin{eqnarray*}}

\newcommand{\beq}{\begin{equation}}
\newcommand{\eeq}{\end{equation}}

\newcommand{\bd}{\begin{enumerate}}
\newcommand{\ed}{\end{enumerate}}

\renewcommand{\>}{\rightarrow}

\usepackage{fancyhdr}
\renewcommand{\leftmark}{{\scriptsize Effective vanishing theorems for vector bundles}}

\begin{document}
\title{Effective vanishing theorems for ample and globally generated vector bundles}
\makeatletter
\let\uppercasenonmath\@gobble
\let\MakeUppercase\relax
\let\scshape\relax
\makeatother

\author{Kefeng Liu}

\address{Center of Mathematical Science, Zhejiang
University, Hangzhou, 310027, CHINA. Department of Mathematics,
UCLA, Los Angeles, CA, USA} \email{liu@cms.zju.edu}

\author{Xiaokui Yang}
\address{Academy of Mathematics and Systems Science, The Chinese Academy of Sciences, Beijing,  CHINA}
\email{xkyang@amss.ac.cn}

\maketitle
\begin{abstract}By proving an integral formula of the curvature tensor of $E\ts \det E$, we observe that the curvature of
$E\ts \det E$ is very similar to that of a line bundle and obtain
certain new Kodaira-Akizuki-Nakano type vanishing theorems for
vector bundles. As special cases, we deduce new vanishing theorems
for ample, nef and globally generated vector bundles by analytic
method instead of the Leray-Borel-Le Potier spectral sequence.
\end{abstract}

\maketitle

\section{Introduction}
Many vanishing theorems have been obtained for the
 Dolbeault cohomology of ample and globally generated vector
 bundles on smooth projective manifolds, mainly due to the efforts
 of J. Le Potier, M. Schneider, T. Peternell,  A.~J.~Sommese, J-P. Demailly, L. Ein and
 R. Lazarsfeld, L. Manivel, F.~Laytimi and W. Nahm(\cite{D1, EL, LN, LN2, LN3, Le, M, PLS, So}).

\btheorem[{\cite[Le Potier]{Le}}]  If $E$ is an ample vector bundle
over a smooth projective manifold $X$, then $$H^{p,q}(X,E)=0
\qtq{for any} p+q\geq n+r$$ where $n=\dim_\C X$ and
$r=\mathrm{rank}(E)$. \etheorem

When $r\leq n$, the vanishing pairs $(p,q)$ are contained in a
\textbf{triangle} enclosed by three lines $p+q=n+r$, $p=n$ and
$q=n$. By using the Leray-Borel-Le Potier spectral sequence, many
interesting generalizations are ob\-tained for products of symmetric
and skew-symmetric powers of an ample vector bundle, twisted by a
suitable power of its determinant line bundle. Here, we just list
some results for the reader's convenience, and for more related
results we refer to \cite{D1, EL, LN, LN2, LN3, Le, M} and the
references therein.

\btheorem[{\cite[Manivel]{M}}]\label{G} Let $E$ be ample and $L$ be
nef, or $E$ be nef and $L$ ample.
 $$H^{p,q}(X, S^k E\ts(\det E)^{n-p+1}\ts L)=0 \qtq{for any}
p+q>n$$ where $r$ is the rank of $E$. \etheorem

%
%
%

 The common feature of their results is that the vanishing theorems hold
for $(p, q)$ lying inside or on certain {triangles}.

 As is well-known, except  Nakano's vanishing theorem,  few vanishing theorems for vector bundles are proved by analytic method. In this paper,  we use analytic method to prove vanishing theorems
for certain Dolbeault cohomology groups of  bounded vector bundles.
The new vanishing theorems have quite different features and they
hold  for $(p, q)$ lying inside or on certain symmetric
\textbf{quadrilaterals}.

\bdefinition\label{bounds}  Let $E$ be an arbitrary holomorphic
vector bundle with rank $r$, $L$
 an ample line bundle and $\eps_1,\eps_2\in \R$. $E$ is said to be
$(\eps_1,\eps_2)$-bounded by $L$ if there exists a Hermitian metric
$h$ on $E$ and a \emph{positive} Hermitian metric $h^L$ on $L$ such
that the curvature $\Theta^{E,h}$ of $E$ is bounded by the
curvatures of $L^{\eps_1}$ and $L^{\eps_2}$, i.e. \beq
\eps_1\omega_L\ts Id_E\leq \Theta^{E,h}\leq \eps_2\omega_L\ts
Id_E\label{bounded} \eeq in the sense of Griffiths.  $E$ is called
\emph{strictly} $(\eps_1,\eps_2)$-bounded by $L$ if, in addition, at
least one of
 $\Theta^{E,h}-\eps_1\omega_L\ts Id_E$ and $\Theta^{E,h}-\eps_2\omega_L\ts Id_E$ is not identically zero. \edefinition

  It is easy to see  that, if $E$ is
  $(\eps_1,\eps_2)$-bounded by $L$, then
  $$E\ts L^{-\eps_1} \qtq{and} E^*\ts L^{\eps_2}$$
are semi-positive in the sense of Griffiths.  In particular, if
$\det E$ is ample, one can choose $L=\det E$ as a natural bound for
$E$. Hence, Definition \ref{bounds} works naturally for many vector
bundles in algebraic geometry. We list some examples as follows. See
Proposition \ref{naturalbounds} for more details.
 \bd \item If $E$ is globally
generated, $E$ is  strictly $(0,1)$-bounded by $L\ts \det E$ for any
ample line bundle $L$;

\item If $E$ is an ample vector bundle with rank $r$,  then $E$ is strictly
$(-1,r)$-bounded by $\det E$;

\item If $E$ is nef with rank $r$,  then $E$ is strictly
$(-1,r)$-bounded by $L\ts \det E$ for arbitrary ample line bundle
$L$;

\item If $E$ is Griffiths-positive,  $E$ is strictly $(0,1)$-bounded
by $\det E$.

\ed

Now we describe our main results  briefly.

\btheorem\label{main1} Let $L$ be ample. If $E$ is strictly
$(\eps_1,\eps_2)$-bounded by $L$ and
 $m+(r+k)\eps_1> 0$, then \beq H^{p,q}(X,S^kE\ts
\det E\ts L^m)=H^{q,p}(X,S^kE\ts \det E\ts L^m)=0, \eeq if $p\geq 1,
q\geq 1$ satisfy \beq \min\left\{ \frac{n-q}{p},
\frac{n-p}{q}\right\}\leq \frac{m+(r+k)\eps_1}{m+(r+k)\eps_2}
.\label{qua7}\eeq In particular, $S^kE\ts \det E\ts L^m$ is both
Nakano-positive and dual-Nakano-positive. Hence
$$ H^{n,q}(X, S^kE\ts \det E\ts L^m)=H^{q,n}(X, S^k{E}\ts \det E\ts
L^m)=0 $$ for $q\geq 1$. \etheorem

 \bremark\bd\item\ $(p,q)$ satisfies condition (\ref{qua7}) if only
if it lies inside or on the following quadrilateral
$Q=A_0A_1A_2A_3$. See Figure $1$ with $A_0$ and $A_2$ removed. Here
$A_0= (0,n)$, $A_1=(n,n)$, $A_2=(n,0)$, $A_3=(c_0,c_0)$ and \beq
c_0= \frac{n}{1+\frac{m+(r+k)\eps_1}{m+(r+k)\eps_2}}.  \eeq  It is
obvious that $Q$ is symmetric with respect to the line $p=q$.
\\
{\setlength{\parindent}{1cm} \setlength{\unitlength}{0.27cm}
\begin{picture}(12,16)(1,0)
{\linethickness{0.015mm} \multiput(0,0)(1,0){15} {\line(0,1){14}}
\multiput(0,0)(0,1){15} {\line(1,0){14}}}
 \linethickness{0.35mm} \put(14,0){\line(-1,2){4.65}}  \linethickness{0.35mm}
 \put(0,14){\line(2,-1){9.35}}
 \linethickness{0.35mm} \put(0,14){\line(1,0){14}}\linethickness{0.35mm} \put(14,0){\line(0,1){14}}
\put(14,14){$A_1$} \put(-0.9,14.3){$A_0$} \put(14,0){$A_2$}
\put(8.1,8.1){$A_3$} \put(5,-2){Figure $1$}

 \quad {\linethickness{0.015mm}
\multiput(20,0)(1,0){15} {\line(0,1){14}} \multiput(20,0)(0,1){15}
{\line(1,0){14}}}
 \linethickness{0.35mm} \put(34,0){\line(-1,2){4.65}}  \linethickness{0.35mm}
 \put(20,14){\line(2,-1){9.35}}
 \linethickness{0.35mm} \put(20,14){\line(1,0){14}}\linethickness{0.35mm} \put(34,0){\line(0,1){14}}
\put(19,14.3){$A_0$} \put(34,14){$A_1$} \put(34,0){$A_2$}
\put(28.1,8.1){$A_3$}

\put(29.4,9.35){\line(1,-1){5.6}}

\put(29.4,9.35){\line(-1,1){5.6}} \put(25,-2){Figure $2$}

\put(35,3.5){$p+q=n+s_0$}
\end{picture}}
\vspace{0.8cm}

\item The condition $m+(r+k)\eps_1>0$ is necessary, which guarantees that
the vector bundle $S^kE\ts \det E\ts L^m$ is  Griffiths-positive. In
fact, in terms of Hermitian metrics,
\begin{align*}
S^kE\ts \det E\ts L^m&=S^k(E\ts L^{-\eps_1})\ts \det(E\ts
L^{-\eps_1})\ts L^{m+(r+k)\eps_1}\\
&\geq L^{m+(r+k)\eps_1}
\end{align*}
and similarly $S^kE\ts \det E\ts L^m\leq L^{m+(r+k)\eps_2}$. On the
other hand, we will see that the bundle $S^kE\ts \det E\ts L^m$ has
a nice metric $h$ such that $(S^kE\ts \det E\ts L^m, h)$ behaves
very similarly to a positive Hermitian ``line bundle" $(\sL,h_0)$.
Moreover, $m+(r+k)\eps_1$ and $m+(r+k)\eps_2$ are the minimal and
maximal eigenvalues of the curvature of $(\sL,h_0)$ respectively.
From these, one can see that Theorem \ref{main1} is optimal.

\item When $\eps_1$ is very close to $\eps_2$, $E$  is approximate Hermitian-Einstein(\cite[Chapter IV, Section~5]{K}) and so it is semi-stable with respect to $L$ (\cite[Chapter V, Theorem ~8.6]{K}). Moreover,
$H^{p,q}(X,S^kE\ts \det E\ts L^m )=0$ for any $p+q\geq n+1$.

\item If $\eps_1\leq 0$, $\eps_2\geq 0$, and $F$ is an arbitrary nef line
bundle, Theorem \ref{main1} also holds for $S^kE\ts \det E\ts L^m\ts
F$. \ed \eremark

 As applications, we obtain

\btheorem\label{main} If $E$ is a globally generated vector bundle
with rank $r$ and $L$ is an ample line bundle, then for any $k\geq
1, m\geq 1$,
$$ H^{p,q}(X, S^kE\ts (\det E)^m\ts L)=H^{q,p}(X, S^k{E}\ts (\det
E)^m\ts L)=0 $$ if  $p\geq 1, q\geq 1$ satisfy \beq
 \min\left\{ \frac{n-q}{p},
\frac{n-p}{q}\right\}\leq \frac{m-1}{m-1+(r+k)}. \label{cha}\eeq In
particular, $S^kE\ts (\det E)^m\ts L$ is both Nakano-positive and
dual-Nakano-positive and
$$ H^{n,q}(X, S^kE\ts (\det E)^m\ts L)=H^{q,n}(X, S^k{E}\ts (\det
E)^m\ts L)=0 $$ for any $q\geq 1$. \etheorem

The right hand side of (\ref{cha}) depends only on the ratios and it
makes Theorem~\ref{main}  quite different from the results in
\cite{D1, M, LN, LN2}. More precisely, for some specific vanishing
pair $(p,q)$, the power of $\det E$ may be independent of the
dimension of $X$. For example, let $n=3n_0+2$, $k=1$ and $m=r+2$. By
(\ref{cha}), we can choose two different pairs $(p,q)=(2,n-1)$ and
$(p,q)=(2n_0+2,2n_0+1)$, and obtain
\begin{align}\label{dif}
&H^{2,n-1}(X,E\ts (\det E)^{r+2}\ts L)\\
=\;&0=H^{2n_0+2,2n_0+1}(X,E\ts(\det E)^{r+2}\ts L) \notag
\end{align}
for any globally generated $E$ and ample $L$.  In general, we do not
have $H^{p,q}(X,E\ts (\det E)^{r+2}\ts L)=0$ \textbf{for all}
$p+q\geq n+1$, if $1<r\ll n$ (One can compare it with Theorem
\ref{G} and also Corollary $1.5$ in \cite{LN2}).
  On the other hand, for fixed
$(k,m)$, the quadrilateral $Q$ contains a triangle $p+q\geq n+s_0$
for some $s_0\in (0, n]$. See Figure $2$. Moreover, if the power $m$
of $\det E$ is large enough, we obtain $H^{p,q}(X,S^kE\ts (\det
E)^m\ts L)=0$ for $p+q\geq n+1$. Examples in \cite{PLS} and
\cite{D1} indicate that a sufficient large power of $\det E$ is
necessary in this case.
 For more details, see Corollary
\ref{1.2}, Corollary \ref{1.3} and Example \ref{D}.

 \btheorem\label{ample} Let $r=rank(E)$. If $E$ is ample (resp. nef)
and $L$ is nef (resp. ample), then for any $k\geq 1$ and $m\geq
k+r+1$,
$$ H^{p,q}(X, S^kE\ts (\det E)^m\ts L)=H^{q,p}(X, S^k{E}\ts (\det
E)^m\ts L)=0,$$ if $p\geq 1, q\geq 1$ satisfy \beq
 \min\left\{ \frac{n-q}{p},
\frac{n-p}{q}\right\} \leq \frac{(m-1)-(r+k)}{(m-1)+r(r+k)}.\eeq
\etheorem

By a similar setting as of (\ref{dif}), it is easy to see that the
result in Theorem~\ref{ample} is different from  the results of
\cite{D1, M, LN, LN2}.

\begin{remark} Our method is a generalization of the
analytic proof of the Kodaira-Akizuki-Nakano vanishing Theorem for
line bundles. We have obtained similar results for certain
``partially" positive vector bundles.
\end{remark}

\noindent\textbf{Acknowledgements.} This paper is based on the
second named author's doctoral thesis at the University of
California at Los Angeles. The authors would like thank the referees
for many suggestions.

\section{Background material}

Let $E$ be a holomorphic vector bundle over a compact K\"ahler
manifold $X$ and $h$ a Hermitian metric on $E$. There exists a
unique connection $\nabla$ which is compatible with the
 metric $h$ and complex structure on $E$. It is called the Chern connection of $(E,h)$. Let $\{z^i\}_{i=1}^n$ be  the local holomorphic coordinates
  on $X$ and  $\{e_\alpha\}_{\alpha=1}^r$ be a local frame
 of $E$. The curvature tensor $R\in \Gamma(X,\Lambda^2T^*X\ts E^*\ts E)$ has the form
 \beq R=\frac{\sq}{2\pi} R_{i\bar j\alpha}^\gamma dz^i\wedge d\bar z^j\ts e^\alpha\ts e_\gamma\eeq
where $R_{i\bar j\alpha}^\gamma=h^{\gamma\bar\beta}R_{i\bar
j\alpha\bar \beta}$ and \beq R_{i\bar j\alpha\bar\beta}= -\frac{\p^2
h_{\alpha\bar \beta}}{\p z^i\p\bar z^j}+h^{\gamma\bar
\delta}\frac{\p h_{\alpha \bar \delta}}{\p z^i}\frac{\p
h_{\gamma\bar\beta}}{\p \bar z^j}.\eeq Here and henceforth we
 adopt the Einstein convention for summation.

\bdefinition
 A Hermitian vector bundle
$(E,h)$ is said to be \emph{Griffiths-positive}, if for any nonzero
vectors $u=u^i\frac{\p}{\p z^i}$ and $v=v^\alpha e_\alpha$,  \beq
\sum_{i,j,\alpha,\beta}R_{i\bar j\alpha\bar \beta}u^i\bar u^j
v^\alpha\bar v^\beta>0.\eeq $(E,h)$ is said to be
\emph{Nakano-positive}, if for any nonzero vector
$u=u^{i\alpha}\frac{\p}{\p z^i}\ts e_\alpha$, \beq
\sum_{i,j,\alpha,\beta}R_{i\bar j\alpha\bar \beta} u^{i\alpha}\bar
u^{j\beta}>0. \eeq $(E,h)$ is said to be
\emph{dual-Nakano-positive}, if for any nonzero vector
$u=u^{i\alpha}\frac{\p}{\p z^i}\ts e_\alpha$, \beq
\sum_{i,j,\alpha,\beta}R_{i\bar j\alpha\bar \beta} u^{i\beta}\bar
u^{j\alpha}>0. \eeq It is easy to see that $(E,h)$ is
dual-Nakano-positive if and only if $(E^*,h^*)$ is Nakano-negative.

The notions of semi-positivity, negativity and semi-negativity can
be defined similarly. We say $E$ is Nakano-positive (resp.
Griffiths-positive, dual-Nakano-positive, $\ldots$), if it admits a
Nakano-positive (resp. Griffiths-positive, dual-Nakano-positive,
$\ldots$) metric. \edefinition

 The following analytic definition of nefness is due to
\cite{DPS}.

\bdefinition
 Let $(X,\omega_0)$ be a compact K\"ahler manifold. A line bundle $L$ over $X$ is said to be nef,
if  for any $\eps>0$, there exists a smooth Hermitian metric
$h_\eps$ on $L$ such that the curvature of $(L,h_\eps)$ satisfies
\beq R=-\frac{\sq}{2\pi}\p\bp\log h_\eps\geq -\eps\omega_0.\eeq
\edefinition

 This means that the curvature of $L$ can have an arbitrarily small negative
part. Clearly a nef line bundle $L$ satisfies $$\int_C c_1(L)\geq
0$$ for all irreducible curves $C\subset X$.
 For projective
algebraic $S$ both notions coincide.

 By the Kodaira embedding theorem, we have the following
analytic definition of ampleness.

\bdefinition Let $(X,\omega_0)$ be a compact K\"ahler manifold. A
line bundle $L$ over $X$ is said to be ample, if there exists a
smooth Hermitian metric $h$ on $L$ such that the curvature $R$ of
$(L,h)$ satisfies \beq R=-\frac{\sq}{2\pi}\p\bp\log h>0.
 \eeq

\edefinition \bdefinition

 Let $E$ be a Hermitian vector bundle of rank $r$ over a compact K\"ahler manifold $X$, $L=\sO_{\P(E^*)}(1)$ be the
  tautological line bundle on the projective bundle $\P(E^*)$ and $\pi$  the canonical projection $\P(E^*)\> X$.
  By definition(\cite{Har}), $E$ is an ample vector bundle over $X$ if
   $\sO_{\P(E^*)}(1)$ is an ample line bundle over $\P(E^*)$. $E$ is said to be
   nef, if $\sO_{\P(E^*)}(1)$ is nef.

\edefinition

\bdefinition Let $E$ be a holomorphic vector bundle over $X$. $E$ is
called globally generated, if there exist global holomorphic
sections $\sigma_1,\ldots, \sigma_N\in H^0(X,E)$ such that for all
$x\in X$, $\sigma_1(x),\ldots, \sigma_N(x)$ generate $E_x$. It is
obvious that every globally generated vector bundle possesses a
smooth Hermitian metric with semi-Griffiths positive
curvature(\cite[Corollary~11.5]{D}) \edefinition

 For  comprehensive descriptions of positivity, nefness,
ampleness and related topics, see \cite{D, DPS, L}. In particular,
one can see Section $2$ and Section~$3$ in \cite{LSYY}.

In the following, we will describe the idea of proving vanishing
theorems by using an analytic method. Let $(\phi_{i\bar j})_{n\times
n}$ be a Hermitian positive matrix with eigenvalues \beq
\lambda_1\leq \cdots \leq \lambda_n. \eeq Let $u=\sum u_{I\bar
J}dz^I\wedge d\bar z^J$ be a $(p,q)$ form on $\C^n$ where $u_{I\bar
J}$ is alternate in the indices $I=(i_1,\ldots, i_p)$ and
$J=(j_1,\ldots, j_q)$. We define \beq
T(u,u)=\la[\phi,\Lambda_\omega]u,u\ra \eeq where $\phi=\sq
\phi_{i\bar j}dz^i\wedge d\bar z^j$ and $\Lambda_\omega$ is the
contraction operator of the standard K\"ahler metric on $\C^n$. The
following linear algebraic result is obvious(\cite{D}, p. $334$):
\blemma We have the following estimate \beq T(u,u)\geq
\max\{p\lambda_1-(n-q)\lambda_n, q\lambda_1- (n-p)\lambda_n\} |u|^2.
\label{1}\eeq \elemma

 The following result is well-known.

\bcorollary \label{kodaira} Let $(L,h)$ be a Hermitian line bundle
over a compact K\"ahler manifold $(X,\omega_0)$. Let $\lambda_1$ and
$\lambda_n$ be the smallest and largest eigenvalue functions of
$R^L$ with respect to $\omega_0$ respectively. If
$$\max\{p\lambda_1-(n-q)\lambda_n, q\lambda_1- (n-p)\lambda_n\}$$ is
positive everywhere, then \beq H^{p,q}(M,L)=H^{q,p}(M,L)=0.
\eeq\bproof  By a well-known Bochner formula for $L$,
$$\Delta''=\Delta'+[R^L,\Lambda_{\omega_0}]$$
for any $u\in \Om^{p,q}(M,L)$, \beq \la \Delta'' u,u\ra=\la\Delta'
u, u\ra+T(u,u). \eeq If $\Delta''u=0$, we get $u=0$ since
$T(u,u)\geq 0$. \eproof \ecorollary

\bremark The condition in Corollary \ref{kodaira} can be satisfied
if and only if $(L,h)$ is Griffiths positive or Griffiths-negative.
If $(L,h)$ is a positive line bundle over a compact complex manifold
$X$, we can define a K\"ahler metric on $X$\beq \omega_0=
R^L=-\frac{\sq}{2\pi}\p\bp\log h. \eeq In this case, $\phi= R^L$ in
Lemma \ref{kodaira} and $\lambda_1=\lambda_n=1$. Hence, if $p+q\geq
n+1$, $H^{p,q}(X,L)=0$. This is the Kodaira-Akizuki-Nakano vanishing
theorem. But in general, if $R^L$ is not related to $\omega_0$, we
can only get a part of vanishing cohomology groups by this method.
More precisely, we can only obtain a vanishing quadrilateral as
Figure $1$. \eremark
 Let $(E,h)$ be a
Hermitian holomorphic vector bundle with rank $r$ over a compact
K\"ahler manifold $(X,\omega_g)$. For any fixed point $p\in X$,
there exist local holomorphic coordinates  $\{z^i\}_{i=1}^n$ and
local holomorphic frames $\{e_\alpha\}_{\alpha=1}^r$ such that \beq
g_{i\bar j}(p)=\delta_{ij},\quad
h_{\alpha\bar\beta}(p)=\delta_{\alpha\bar\beta}.\eeq
The curvature term in the formula $\Delta''=\Delta'+[R^E,\Lambda_g]$
can be written as
\begin{align}\label{curvatureformula}
T(u,u)&=\la  [R^E, \Lambda_g]u,u\ra\\
\notag &= \sum_{}R_{i\bar j \alpha\bar\beta}u_{I,\bar
{iS},\alpha}\bar u_{I,\bar {jS}\beta}+\sum_{} R_{i\bar
j\alpha\bar\beta} u_{jR,\bar J,\alpha}\bar u_{iR,\bar
J,\beta}\\
\notag &\quad -\sum_{}R_{ii\alpha\bar\beta}u_{I\bar J\alpha}\bar
u_{I\bar J\beta} \end{align} for any $u=\sum u_{I\bar J\alpha
}dz^I\wedge d\bar z^J\ts e_\alpha$. For more details, see (\cite{D},
p. $341$). From formula (\ref{curvatureformula}), it is very
difficult to obtain vanishing theorems for vector bundles. If the
curvature $R^E$ has a nice expression, for example \beq R_{i\bar
j\alpha\bar\beta}=\phi_{i\bar
j}\tau_\alpha\bar\tau_\beta\label{positivecurvature} \eeq then $E$
behaviors as a line bundle with curvature $(\phi_{i\bar j})$.
Unfortunately, few examples with property (\ref{positivecurvature})
can be found( Note also that the curvature formulation here is
stronger than the curvature of projectively flat vector bundles).
However, an integral version of (\ref{positivecurvature}) exists on
vector bundles of type $E\ts \det E$,
\begin{align}\label{curv2}
R^{E\ts \det E}_{i\bar j\alpha\bar\beta}(s)&=R_{i\bar j
\alpha\bar\beta}(s)+\delta_{\alpha\beta}\cdot \sum_{\gamma}R_{i\bar
j \gamma\bar\gamma}(s)\\
\notag &=r! \cdot \int_{ \P^{r-1}} \frac{\phi_{i\bar j} W_\alpha\bar
W_\beta}{|W|^{2}} \frac{\omega^{r-1}_{FS}}{(r-1)!}
\end{align}
where $[W_1,\ldots,W_r]$ are the homogeneous coordinates on
$\P^{r-1}$, $\omega_{FS}$ is the Fubini-Study metric and \beq
\phi_{i\bar j}=(r+1)\sum_{\gamma,\delta}R_{i\bar j
\gamma\bar\delta}(s)\frac{W_\delta\bar W_\gamma}{|W|^2}. \eeq It is
obvious that if $E$ is Griffiths-positive, then $E\ts \det E$ is
both Nakano-positive and dual-Nakano-positive. With the help of the
nice formulation (\ref{curv2}), we  obtain vanishing theorems
similar to Corollary \ref{kodaira} for vector bundles.

\section{Vanishing theorems for bounded vector bundles}

Firstly, we would like to recall the following \bdefinition

Let $E$ be an arbitrary holomorphic vector bundle with rank $r$, $L$
 an ample line bundle and $\eps_1,\eps_2\in \R$. $E$ is said to be
$(\eps_1,\eps_2)$-bounded by $L$ if there exists a Hermitian metric
$h$ on $E$ and a \emph{positive} Hermitian metric $h^L$ on $L$ such
that the curvature of $E$ is bounded by the curvatures of
$L^{\eps_1}$ and $L^{\eps_2}$, i.e. \beq \eps_1\omega_L\ts Id_E\leq
\Theta^{E,h}\leq \eps_2\omega_L\ts Id_E\label{bounded1} \eeq in the
sense of Griffiths.  $E$ is called \emph{strictly}
$(\eps_1,\eps_2)$-bounded by $L$ if, in addition, at least one of
 $\Theta^{E,h}-\eps_1\omega_L\ts Id_E$ and $\Theta^{E,h}-\eps_2\omega_L\ts Id_E$ is not identically zero. \edefinition

 It is easy to see that $E$ is $(\eps_1,\eps_2)$-bounded by
$L$ if and only if $E\ts L^{-\eps_1}$
 and $E^*\ts L^{\eps_2}$ are Griffiths-semi-positive. Similarly, if $E$ is strictly
$(\eps_1,\eps_2)$-bounded by $L$, then at least one of the
Griffiths-semi-positive vector bundles $E\ts L^{-\eps_1}$
 and $E^*\ts L^{\eps_2}$ is  not trivial.

\bproposition\label{naturalbounds} Let $E$ be a holomorphic vector
bundle  with rank $r$ over a projective manifold.\bd
\item If $E$ is globally generated, $E$ is  strictly $(0,1)$-bounded by $L\ts \det E$ for any ample line
bundle $L$;

\item If $E$ ample,  $E$ is strictly
$(-1,r)$-bounded by $\det E$;

\item If $E$ is nef, $E$ is strictly
$(-1,r)$-bounded by $L\ts \det E$ for any ample line bundle $L$;

\item If $E$ is Griffiths-positive, $E$ is strictly $(0,1)$-bounded
by $\det E$.

\ed \bproof $(1)$ As is well-known, if $E$ is globally generated,
there exists a Hermitian metric $h$ on $E$ such that $\Theta^{E,h}$
is Griffiths-semi-positive and  $E\ts \det E^*=\Lambda^{r-1}E^*$ is
Griffiths-semi-negative.  If $L$ is an ample line bundle,  $E\ts
\det E^*\ts L^*$ is Griffiths-negative, i.e.
$$\Theta^{E,h}<\omega_{L\ts\det E}\ts Id_E.$$
Hence, $E$ is strictly $(0,1)$-bounded by $L\ts \det E$.

$(2)$ We assume $r>1$. By a result of \cite{Bo1}, \cite{MT1} and
\cite{LSY13}, if $E$ is ample, $E\ts \det E$ is Griffiths-positive.
On the other hand, $E^*\ts \det E=\Lambda^{r-1}E$ is ample and so
$(E^*\ts\det E)\ts \det(E^*\ts \det E)=E^*\ts (\det E)^r$ is
Griffiths-positive.

$(3)$ If $E$ is nef, $S^{r+1}E\ts L$ is ample and by a result of
\cite{LSY13}, $E\ts \det E\ts L$ is Griffiths-positive. Similarly,
we know $S^{r+1}(E^*\ts \det E)\ts L$ is ample and so $E^*\ts (\det
E)^r\ts L$ is Griffiths-positive.

$(4)$ It is obvious. \eproof \eproposition

\bremark In general, if $E$ is $(-1,r)$-bounded by $\det E$, $E$ is
not necessarily ample. For example, let $E=L^3\ds L^{-1}$ for some
ample line bundle $L$, then $E$ is $(-1,2)$ bounded by $\det E=L^2$.

\eremark

 Let $\omega_{FS}$ be the standard Fubini-Study metric on $\P^{r-1}$ with $\int_{\P^{r-1}}\omega_{FS}^{r-1}=1$ and
 $[W_1,\ldots W_r]$ the homogeneous coordinates on $\P^{r-1}$.
  If $A=(\alpha_1,\ldots, \alpha_k)$
and $B=(\beta_1,\beta_2,\ldots, \beta_k)$, we  define the
generalized Kronecker-$\delta$ for multi-index by the following
formula \beq \delta_{AB}=\sum_{\sigma\in S_k}
\prod_{j=1}^k\delta_{\alpha_{\sigma(j)}\beta_{\sigma(j)}}
\label{generalizeddelta}\eeq where $S_k$ is the permutation group in
$k$ symbols. The following linear algebraic lemma is obvious (see
also \cite{LSY13}). \blemma\label{calculation} Let $W=[W_1,\ldots,
W_r]$ be the homogeneous coordinates on $\P^{r-1}$. If
$V_A=W_{\alpha_1}\cdots W_{\alpha_k}$and $V_B=W_{\beta_1}\cdots W_{
\beta_k}$, then \beq\int_{\P^{r-1}} \frac{V_{A}\bar
V_B}{|W|^{2k}}\frac{\omega_{FS}^{r-1}}{(r-1)!}=\frac{\delta_{AB}}{(r+k-1)!}.
\eeq For simple-index notations, \beq\begin{aligned} \int_{\P^{r-1}}
\frac{W_\alpha \bar
W_\beta}{|W|^2}\frac{\omega_{FS}^{r-1}}{(r-1)!}&=\frac{\delta_{\alpha\beta}}{r!}
,\\[1ex]
\int_{\P^{r-1}}\frac{W_\alpha \bar{W_\beta} W_\gamma
\bar{W_\delta}}{|W|^4}\frac{\omega_{FS}^{r-1}}{(r-1)!}&=\frac{\delta_{\alpha\beta}\delta_{\gamma\delta}+\delta_{\alpha\delta}\delta_{\beta\gamma}}{(r+1)!}.
\end{aligned}
\eeq

%
\elemma

Let $h$ be a Hermitian metric on the vector bundle $E$. At a fixed
point $p\in X$, if we assume
$h_{\alpha\bar\beta}=\delta_{\alpha\bar\beta}$, then the naturally
induced bundle $(E\ts (\det E)^m, h\ts (\det h)^m)$ has curvature
components \beq R_{i\bar j \alpha\bar\beta}^{E\ts (\det
E)^m}=R_{i\bar j \alpha\bar\beta}+\delta_{\alpha\beta}\cdot
m\sum_\delta R_{i\bar j \delta\bar\delta} \eeq where $R_{i\bar j
\alpha\bar\beta}$ is the curvature component of $(E,h)$. It is
obvious that $S^kE$ has basis \beq \{e_A= e_{1}^{\alpha_1}\ts \cdots
\ts e^{\alpha_r}_r\} \eeq if $A=(\alpha_1,\ldots, \alpha_r)$  with
$\alpha_1+\cdots+\alpha_r=k$ and $\alpha_j$  are nonnegative
integers. Similarly, $(S^kE\ts (\det E)^m, S^kh\ts (\det h)^m)$ has
curvature components \beq R_{i\bar j A\bar B}^{S^kE\ts (\det
E)^m}=R_{i\bar j A\bar B}+\delta_{AB}\cdot m\sum_\delta R_{i\bar j
\delta\bar\delta} .\eeq

  \blemma\label{linear}\hspace{-1ex} If $(E,h)$ is a Hermitian vector bundle,  the
curvature of $(S^kE\ts (\det E)^m, S^kh\ts (\det h)^m)$ can be
written as\beq R^{S^kE\ts (\det E)^m}_{i\bar j A\bar B}(p)=(r+k-1)!
\cdot \int_{ \P^{r-1}} \frac{ V_A\bar V_B}{|W|^{2k}} \phi_{i\bar
j}\frac{\omega^{r-1}_{FS}}{(r-1)!} \eeq where \beq \phi_{i\bar
j}=(r+k)\sum_{\gamma,\delta}R_{i\bar j
\gamma\bar\delta}(p)\frac{W_\delta\bar
W_\gamma}{|W|^2}+(m-1)\sum_\delta R_{i\bar j
\delta\bar\delta}.\label{phi} \eeq \bproof It follows from Lemma
\ref{calculation}.\eproof \elemma 

Now we prove Theorem \ref{main1}. \btheorem\label{maintheorem}  Let
$L$ be ample. If $E$ is strictly $(\eps_1,\eps_2)$-bounded by $L$
and
 $m+(r+k)\eps_1> 0$, then \beq H^{p,q}(X,S^kE\ts
\det E\ts L^m)=H^{q,p}(X,S^kE\ts \det E\ts L^m)=0 \eeq if $p\geq 1,
q\geq 1$ satisfy \beq \min\left\{ \frac{n-q}{p},
\frac{n-p}{q}\right\}\leq \frac{m+(r+k)\eps_1}{m+(r+k)\eps_2}
.\label{qua}\eeq In particular, if $m\!+\!(r\!+\!k)\eps_1\!>\!0$,
$S^kE\!\ts\! \det E\!\ts\! L^m$ is both Nakano-positive and
dual-Nakano-positive  and
$$ H^{n,q}(X, S^kE\ts \det E\ts L^m)=H^{q,n}(X, S^k{E}\ts \det E\ts
L^m)=0 $$ for $q\geq 1$. \bproof Let $h$ be a Hermitian metric on
$E$ and $h^L$ a positive Hermitian metric on $L$ such that
$$\eps_1\omega_L\ts Id_E\leq \Theta^{E,h}\leq \eps_2\omega_L\ts
Id_E.$$ We can polarize $X$ by \beq \omega_g=\omega_L=-\frac{\sq
}{2\pi}\p\bp\log  h^L.\eeq
 At a fixed point $z_0\in X$, we can assume
$$ g_{i\bar j}(z_0)=\delta_{i\bar j}   \qtq{and} h_{\alpha\bar
\beta}(z_0)=\delta_{\alpha\bar\beta}.$$ Therefore, \beq g_{i\bar
j}(z_0)=R_{i\bar j}^{h^L}(z_0)=\delta_{i\bar j}.\label{metric} \eeq
Let \beq \phi_{i\bar j}=(r+k)\left(\sum_{\gamma,\delta}R^h_{i\bar j
\gamma\bar\delta}(z_0)\frac{W_\delta\bar
W_\gamma}{|W|^2}\right)+mR^{h_L}_{i\bar j}\label{phi1} ,\eeq then by
Lemma \ref{linear}, the curvature tensor of $ S^kE\ts \det  E\ts
L^m$ can be written as \beq R^{S^kE\ts \det E\ts L^m}_{i\bar j A\bar
B}(z_0)=(r+k-1)! \cdot \int_{ \P^{r-1}} \frac{ V_A\bar
V_B}{|W|^{2k}} \phi_{i\bar j}\frac{\omega^{r-1}_{FS}}{(r-1)!}
.\label{curv}\eeq
 It is easy
to see that, if $E$ is $(\eps_1,\eps_2)$-bounded by $L$, then for
any $v=(v^1,\ldots, v^n)\in \C^n\setminus\{0\}$, at point $z_0$,
$$\eps_1 R_{i\bar j}^{h_L} v^i\bar{v}^j\leq \sum_{\gamma,\delta}R^h_{i\bar j
\gamma\bar\delta}(z_0)\frac{W_\delta\bar
W_\gamma}{|W|^2}v^i\bar{v}^j\leq \eps_2 R_{i\bar j}^{h_L}
v^i\bar{v}^j. $$ By formula (\ref{phi1}), we obtain
 \beq \left(m+(r+k)\eps_1\right)|v|^2 \leq
\phi_{i\bar j}v^i \bar v^j\leq \left(m+(r+k)\eps_2\right)|v|^2.
\label{ww}\eeq
 Since
$m\!+\!(r\!+\!k)\eps_1\!>\!0$, it is obvious that $S^kE\!\ts\! \det
E\!\ts\! L^m$ is both Nakano-positive and dual-Nakano-positive by
(\ref{curv}).
  Let $\lambda_1$ be the smallest
eigenvalue of $(\phi_{i\bar j})$ and $\lambda_n$ the largest one,
then \beq m+(r+k)\eps_1\leq  \lambda_1\leq \lambda_n\leq
m+(r+k)\eps_2. \label{eigenratio}\eeq Let
$\phi=\frac{\sq}{2\pi}\phi_{i\bar j}dz^i\wedge d\bar z^j$. We
consider the curvature term in the Bochner formula
$$\Delta''=\Delta'+[R,\Lambda_g]$$ for vector bundle  $S^kE\ts \det E\ts L^m$.  For any nonzero
$$u=u_{I\bar JA}dz^I\wedge d\bar z^J\ts e_A\in\Om^{p,q}(X, S^kE\ts
\det E\ts L^m),$$ we set $$U=\sum_A u_{I\bar JA}V_A dz^I\wedge d\bar
z^J,$$ then the curvature term can be written as
 \begin{align*} T(u,u)&=\la [R,\Lambda_g]u,u\ra
 = (r+k-1)!\int_{\P^{r-1}}\left\la
[\phi,\Lambda_{g}]U,U \right\ra\cdot
\frac{1}{|W|^{2k}}\cdot\frac{\omega^{r-1}_{FS}}{(r-1)!}\\
&\geq (r+k-1)!\\
&\quad \cdot\!\!\int_{\P^{r-1}}\!\!\max\{p\lambda_1-(n-q)\lambda_n,
q\lambda_1-(n-p)\lambda_n \}|U|^2\cdot
\frac{1}{|W|^{2k}}\cdot\frac{\omega^{r-1}_{FS}}{(r-1)!}\\
&=\max\{pK_1-(n-q)K_n, qK_1-(n-p)K_n \}
\end{align*}  where
$$ K_i= (r+k-1)! \cdot
\int_{ \P^{r-1}} \frac{ |U|^2}{|W|^{2k}}
\lambda_i\frac{\omega^{r-1}_{FS}}{(r-1)!},\quad  i=1,n.
$$ By (\ref{eigenratio}), if $m+(r+k)\eps_1>0$, then $m+(r+k)\eps_2>0$ and \beq
0< \frac{m+(r+k)\eps_1}{m+(r+k)\eps_2}\leq \frac{K_1}{K_n}.
\label{strict1}\eeq
 Moreover, if $E$ is \emph{strictly} $(\eps_1,\eps_2)$-bounded by $L$,
 we obtain the strict inequality:
\beq 0< \frac{m+(r+k)\eps_1}{m+(r+k)\eps_2}<\frac{K_1}{K_n}
\label{strict}\eeq
 at some point $z_0\in X$.
Hence, if $p\geq 1, q\geq 1$ satisfy \beq \min\left\{ \frac{n-q}{p},
\frac{n-p}{q}\right\}\leq \frac{m+(r+k)\eps_1}{m+(r+k)\eps_2}
\label{3} \eeq we obtain \beq
 \min\left\{ \frac{n-q}{p}, \frac{n-p}{q}\right\}<\frac{K_1}{K_n}.
\eeq By standard Bochner formulas(e.g. Corollary \ref{kodaira}),  we
deduce that $$H^{p,q}(X,S^kE\ts (\det E)^m\ts L)=H^{p,q}(X,S^kE\ts
(\det E)^m\ts L)=0$$ if $(p,q)$ satisfies (\ref{3}). The proof of
Theorem \ref{maintheorem} is complete.\eproof \etheorem

\begin{proof}[Proof of Theorem \ref{main}] When $m=1$, the conclusion is
obvious. In this case,  the only vanishing pair is $(p,q)=(n,n)$ and
it follows from the fact that $H^{n,n}(X,\sE)=0$ if $\sE$ is
Griffiths-positive. Now we consider $m\geq 2$. If $E$ is a globally
generated vector bundle with rank $r$ and $L$ is an ample line
bundle, by Proposition \ref{naturalbounds}, $E$ is strictly
$(0,1)$-bounded by $L^{\frac{1}{m-1}}\ts \det E$ (in the curvature
sense, here we do not really use the concept of $\Q$-line bundles).
Theorem \ref{main} follows from Theorem \ref{main1} for
$(\eps_1,\eps_2)=(0,1)$ and the relation
$$S^kE\ts (\det E)^m\ts L=S^kE\ts \det E\ts \left(L^{\frac{1}{m-1}}\ts \det E\right)^{m-1}.\vspace{-2em}$$
\end{proof}

\begin{proof}[Proof of Theorem \ref{ample}] Similarly, assume $m\geq 2$. If
$E$ is ample (resp. nef) and $L$ is nef (resp. ample), by
Proposition \ref{naturalbounds}, $E$ is strictly $(-1,r)$-bounded by
$L^{\frac{1}{m-1}}\ts \det E$. Theorem \ref{ample} follows from
Theorem \ref{main1} for  $(\eps_1,\eps_2)=\linebreak (-1,r)$.
\end{proof}

Similarly, we have, \btheorem Let $(E,h)$ be a Hermitian vector
bundle with semi-Griffiths positive (resp. Griffiths positive)
curvature and $L$ is an ample (resp. nef) line bundle, we have
 $$ H^{p,q}(X, S^k E\ts (\det  E)^m\ts L)=H^{q,p}(X,
S^k{ E}\ts (\det  E)^m\ts L)=0 $$ if $p\geq 1, q\geq 1$ satisfy
 $$ \min\left\{ \frac{n-q}{p},
\frac{n-p}{q}\right\}\leq \frac{m-1}{r+k+m-1}. $$ \bproof It follows
from part (4) of Proposition \ref{naturalbounds} and Theorem
\ref{main1}. \eproof\etheorem

 Now we want to analyze the condition \beq
\min\left\{ \frac{n-q}{p}, \frac{n-p}{q}\right\}\leq
\lambda_0\label{inequality} \eeq for some $\lambda_0\!\in\! [0,1]$.
Without loss of generality, we assume $p\!\geq\! q\!\geq\!
 1$, then~(\ref{inequality}) is equivalent to \beq p+\lambda_0q\geq n .\eeq When $p=q$, we obtain
\beq c_0=\frac{n}{1+\lambda_0}.\eeq $  (p,q)$ satisfies
(\ref{inequality}) if and only if $(p,q)$ lies in the quadrilateral
$Q\!=\!A_0A_1A_2A_3$ where \beq A_0= (0,n),\quad A_1=(n,n),\quad
A_2=(n,0),\quad A_3=(c_0,c_0). \eeq

In the following, we consider the vanishing triangle  shown in
Figure $2$.

 \bcorollary\label{1.2}Let $E$ be globally generated and $L$ be
ample.
 \bd \item If the pair $(k,m,s)$ satisfies \beq m\geq
\frac{1}{s}\left[\frac{n-s}{2}\right](r+k)+1\eeq where $[\bullet]$
is the integer part of $\bullet$, then
$$ H^{p,q}(X,S^kE\ts (\det E)^m\ts L)=0 $$ for any $p+q\geq n+s$.

\item  For fixed $(k,m)$, we have  $$
H^{p,q}(X,S^kE\ts (\det E)^m\ts L)=0 $$ for any $(p,q)$ satisfies
 \beq p+q\geq n+\left(
\frac{2n}{1+\frac{m-1}{r+k+m-1}}-n\right).\eeq \ed \bproof If
$m\geq\frac{1}{s} \left[\frac{n-s}{2}\right](r+k)+1$, we get \beq
\frac{m-1}{r+k+m-1}\geq
\frac{\left[\frac{n-s}{2}\right]}{\left[\frac{n-s}{2}\right]+s}
.\eeq If $p+q\geq n+s$, \beq\max_{p+q\geq n+s}\min\left\{
\frac{n-q}{p}, \frac{n-p}{q}\right\}=
\frac{\left[\frac{n-s}{2}\right]}{\left[\frac{n-s}{2}\right]+s}.\eeq
 Part $(1)$  follows from Theorem \ref{main}. For part $(2)$, if
$$p+q\geq n+\left(
\frac{2n}{1+\frac{m-1}{r+k+m-1}}-n\right)=\frac{2n}{1+\frac{m-1}{r+k+m-1}}$$
then \beq \max\{p,q\}\geq \frac{n}{1+\frac{m-1}{r+k+m-1}} .\eeq That
is
$$\frac{m-1}{r+k+m-1}\geq \min\left\{ \frac{n-q}{p},
\frac{n-p}{q}\right\}.$$  Hence part (2) follows. \eproof

\ecorollary

\bremark  Theorem \ref{main} and Corollary \ref{1.2} are also valid
for semi-Griffiths positive $E$. Consider the example $E=T\P^2\ts
\sO_{\P^{2}}(-1)$ with the canonical metric. Since $r=n=2$, by
Corollary \ref{1.2}, we obtain  \beq H^{p,q}(X,E\ts (\det E)^m\ts
L)=0 \eeq for any $p+q\geq n+1$ if $m\geq 1$. It is obvious that the
lower bound $1$ is sharp since \beq H^{n,n-1}(X,E\ts L)\cong
H^{1,1}(\P^n,\C)=\C \eeq if we choose $L=\sO_{\P^n}(1)$ and $m=1$.
So the lower bound
$$\frac{1}{s}\left[\frac{n-s}{2}\right](r+k)+1$$ can not be improved
by a universal constant, i.e., a constant independent on $r,s,n,k$.
Hence the lower bound is optimal in that sense.\eremark

 Similarly, we obtain \bcorollary\label{1.3} Let $E$ be
ample (resp. nef) and $L$  be nef (resp. ample). Suppose $k\geq 1$
and $m\geq r+k+1$. \bd\item If the pair $(k,m,s)$ satisfies
$$m\geq \frac{1}{s}\left[\frac{n-s}{2}\right](r+k)(r+1)+(r+1)+k,$$
then $$ H^{p,q}(X,S^kE\ts (\det E)^m)=0 $$ for any $p+q\geq n+s$.

\item  For fixed $(k,m)$, we have $$ H^{p,q}(X,S^kE\ts (\det E)^m\ts L)=0 $$
 for any $(p,q)$ satisfies $$ p+q\geq n+\left(
\frac{2n}{1+\frac{(m-1)-(r+k)}{(m-1)+r(r+k)}}-n\right).$$
  \ed

\ecorollary

\section{Examples}

It is well-known that globally generated vector bundles are
Griffiths semi-positive. On the other hand, any globally generated
vector bundle has a quotient metric induced from the trivial vector
bundle and so it is semi-dual-Nakano-positive(\cite{D}).

\bcorollary\label{GL} Let $E$ be a globally generated vector bundle
and $L$ an ample line bundle over a projective manifold $X$, then
$S^kE\ts L$ is dual-Nakano-positive for any $k\geq 1$. Moreover,
\beq H^{p,n}(X,S^kE\ts L)=0 \eeq for any $p\geq 1$. \ecorollary

 However, in general, we can not obtain a vanishing
quadrilateral for $S^kE\ts L$ as Figure $1$. It is easy to see that
the result in Corollary \ref{GL} is a vertical line on the boundary
of the quadrilateral in Figure $1$. In \cite{PLS}, the authors found
more vanishing elements close to that vertical line. More precisely,
they proved that \beq H^{p,n-1}(X,S^kE\ts L)=0,\qtq{for any} p\geq
r+1. \eeq
 But in general, they proved that  there exists some $1\leq q\leq n$ such that $H^{n,q}(X, S^kE\ts L)\neq 0$. In particular, $S^kE\ts L$
is not necessarily Nakano-positive. For example, $E=T\P^n\ts
\sO_{\P^n}(-1)$ and $L=\sO_{\P^n}(1)$. It is obvious $E$ is globally
generated. When $n\geq 2$, $E\ts L=T\P^n$ is dual-Nakano-positive
but not Nakano-positive. More generally, we have

\bexample[Demailly, \cite{D1}]\label{D} Let  $X= G(r,V)$ be the
Grassmannian of subspaces of codimension $r$ of a vector space $V$,
$\dim_{\C} V= d$, and $E$ the tautological quotient vector bundle of
rank $r$ over $X$. Then $E$ is globally generated  and $L:= \det E$
is very ample. \beq H^{n,q}(X,S^kE\ts \det E)=\begin{cases}0, &q\neq
(r-1)(d-r);\\S^{k+r-d}V\ts \det V, &q=(r-1)(d-r).\end{cases}\eeq
where $n=\dim_\C X=r(d-r)$. If $r=d-1$, then $X=\P^n=\P^{d-1}$ and
$E=T\P^n\ts \sO_{\P^n}(-1)$, $\det E=\sO_{\P^n}(1)$. That is \beq
H^{n,q}(\P^n,S^kT\P^n\ts
\sO_{\P^n}(1-k))=\begin{cases}0,& q\neq n-1;\\
S^{k-1}V\ts \det V, & q=n-1.\end{cases}\eeq Therefore, if $n\geq 2$,
$S^kT\P^n\ts \sO_{\P^n}(1-k)$ can not be Nakano-positive by the
non-vanishing. However, we shall see that for any $\ell\geq 2-k$,
$S^{k}T\P^n\ts \sO_{\P^n}(\ell)$  is both Nakano-positive and
dual-Nakano-positive. Moreover, we  can obtain more vanishing
results about it.
 \eexample

 Let $h_{FS}$
be the Fubini-Study metric on $\P^n$ and it also induces a metric on
$L=\sO_{\P^n}(1)$. It is easy to see that \beq \omega_L\ts Id\leq
\Theta^{T\P^n}\leq 2 \omega_L\ts Id.\eeq So $T\P^n$ is strictly
$(1,2)$-bounded by $L$. Similarly, $H=T\P^n\ts \sO_{\P^n}(-1)$ is
strictly $(0,1)$-bounded by $L$.

\bproposition\label{ex} If $\ell\geq 2-k$,
$S^kT\P^n\ts\sO_{\P^n}(\ell)$ is Nakano-positive and
dual-Nakano-positive and \beq H^{p,q}(\P^n,S^kT\P^n\ts
\sO_{\P^n}(\ell))=H^{q,p}(\P^n,S^kT\P^n\ts \sO_{\P^n}(\ell))=0 \eeq
for any $p\geq 1, q\geq 1$ satisfy \beq
 \min\left\{\frac{n-p}{q},
\frac{n-q}{p}\right\}\leq \frac{\ell+k-1}{\ell+n+2k-1} .\eeq \bproof
Since $\det (T\P^n)=\sO_{\P^n}(n+1)$, we see
$\det(H)=\sO_{\P^n}(1)$. It follows from the relation \beq S^kH\ts
\det H\ts \sO_{\P^n}(\ell+k-1)=S^kT\P^n\ts \sO_{\P^n}(\ell) \eeq and
Theorem \ref{main1} with $(\eps_1,\eps_2)=(0,1)$ and $m=r+k-1$. Here
$\ell+k-1\geq 1$, i.e., $\ell\geq 2-k$ is necessary and optimal by
Example \ref{D}. \eproof \eproposition

\bremark Although $T\P^n$ is not Nakano-positive when $n\geq 2$,
$S^kT\P^n$ is both Nakano-positive and dual-Nakano-positive for any
$k\geq 2$ (see \cite{LSY13}). It is also easy to see that similar
results as Proposition \ref{ex} hold on general flag manifolds.

\eremark


\begin{thebibliography}{99}



\bibitem{Bo1} B.~Berndtsson,
\newblock {\em Curvature of vector bundles associated to holomorphic fibrations}. Ann. of Math. (2), \textbf{169} (2009), no. 2, 531--560.

\bibitem{D} J-P. Demailly,
\newblock {\em Complex analytic and algebraic geometry}.
\newblock Book online \url{http://www-fourier.ujf-grenoble.fr/~demailly/books.html}.

\bibitem{D0} J-P. Demailly,   {\em Vanishing theorems for tensor powers of a positive vector bundle}.
In: Proceedings of the Conference Geometry and Analysis on Manifolds
(edited by T.~Sunada), Lecture Notes in Math., \textbf{1339},
Springer-Verlag, (1988).

\bibitem{D1} J-P. Demailly,  {\em Vanishing theorems for tensor powers of an ample vector bundle}. Invent. Math., \textbf{91} (1988), 203--220.


\bibitem{D3} J-P. Demailly,
{\em Regularization of closed positive currents and intersection
theory}. J.~Algebraic Geom., \textbf{1} (1992), no.~3, 361--409.


\bibitem{DPS} J-P. Demailly, T.~Peternell and M.~Sehneider, {\em Compact complex manifolds with numerically effective tangent bundles}. J.~Alg. Geom., \textbf{3} (1994), 295--345.


\bibitem{EL} L.~Ein and  R.~Lazarsfeld,
\newblock  {\em Syzygies and Koszul cohomology of smooth projective varieties of arbitrary dimension}. Invent. Math., \textbf{111} (1993), no.~1,
51--67.


\bibitem{G} P.~A.~Griffiths,  {\em Hermitian differential geometry, Chern classes and positive vector bundles}. Global Analysis (papers in honor of K. Kodaira), Princeton Univ.
Press, Princeton, (1969), 181--251.


\bibitem{Har} R.~Hartshorne,  {\em Ample vector bundles}. Publ. Math. I.H.E.S., \textbf{29} (1966), 319--350.



\bibitem{K} S.~Kobayashi, {\em Differential geometry of complex vector bundles}.  Princeton University Press, (1987).

\bibitem{LN}  F.~Laytimi and W.~Nahm,
\newblock {\em A generalization of Le Potier's vanishing theorem}.
Manuscripta math., \textbf{113} (2004), 165--189.

\bibitem{LN2} F.~Laytimi and W.~Nahm, {\em On a vanishing problem of Demailly}. Int.
Math. Res. Not.,  \textbf{47} (2005), 2877--2889.

\bibitem{LN3} F.~Laytimi and W.~Nahm, {\em A vanishing theorem}. Nagoya Math.~J.,
\textbf{180} (2005), 35--43.

\bibitem{L} R.~Lazarsfeld,
\newblock {\em Positivity in algebraic geometry I, II}.
\newblock  Ergebnisse der
Mathematik und ihrer Grenzgebiete. 3. Folge\,/\,A Series of Modern
Surveys in Mathematics,  Springer-Verlag, Berlin, (2004).


\bibitem{Le}  J.~Le Potier,  {\em Annulation de la cohomolgie \`a valeurs dans un fibr\'e vectoriel holomorphe positif de rang quelconque}.  Math. Ann., \textbf{218} (1975), no.~1, 35--53.


\bibitem{LSY13} K.-F.~Liu,  X.-F.~Sun and X.-K.~Yang, {\em Positivity and vanishing
theorems for ample vector bundles}. J.~Algebraic Geom., {\bf 22}
(2013), 303--331.

\bibitem{LSYY} K.-F.~Liu,  X.-F.~Sun, X.-K.~Yang and S.-T.~Yau,  {\em Curvatures of
moduli space of curves and applications}.
\href{http://arxiv.org/abs/1312.6932}{arXiv:1312.6932}


\bibitem{LSY2} K.-F.~Liu and X.-K.~Yang,  {\em Curvatures of direct image
sheaves of vector bundles and applications I}. J.~Differential
Geom., {\bf 98} (2014), 117--145.


\bibitem{M} L.~Manivel,  {\em Vanishing theorems for ample vector bundles}. Invent. math., \textbf{127} (1997), 401--416.


\bibitem{MT1} C.~Mourougane and S.~Takayama,  {\em Hodge metrics and positivity of direct images}.
J.~Reine Angew. Math., \textbf{606} (2007), 167--178.


\bibitem{PLS} T.~Peternell and J.~Le Potier and M.~Schneider,  {\em Vanishing theorems,
linear and quadratic normality}. Invent. Math., \textbf{87} (1987),
573--586.


\bibitem{So} A.~J.~Sommese,
\newblock {\em Submanifolds of abelian varieties}.
\newblock Math. Ann., \textbf{233} (1978),
229--256.

\bibitem{Yang} Q.~Yang,
\newblock {\em $(k,s)$-positivity and vanishing theorems for compact K\"ahler
manifolds}. Internat. J.~Math., \textbf{22} (2011), no.~4, 545--576.


\end{thebibliography}
\end{document}